\documentclass[journal]{IEEEtran}
\usepackage{graphics}
\usepackage{graphicx}
\usepackage{epsfig,amssymb}
\usepackage{amsmath}
\usepackage{mathrsfs}
\usepackage{color}
\usepackage{array}
\usepackage{amsfonts,booktabs}
\usepackage{lipsum,amsmath,multicol}
\usepackage{graphicx}
\usepackage{subfigure}
\usepackage{times}
\usepackage{pdfsync}
\usepackage{pgfplots}
\usepackage{pgfplotstable}
\usepackage[subnum]{cases}
\usepackage{filecontents}
\usepackage{multirow}
\usepackage{algorithm}
\usepackage{algpseudocode}
\usepackage{float}
\usepackage{setspace}
\newtheorem{theorem}{Theorem}
\newtheorem{lemma}{Lemma}
\newtheorem{definition}{Definition}

\newcommand{\beq}{\begin{equation}}
\newcommand{\eeq}{\end{equation}}
\newcommand{\beqa}{\begin{eqnarray}}
\newcommand{\eeqa}{\end{eqnarray}}

\newcommand{\paren}[1]{\left(#1\right)}
\newcommand{\sqparen}[1]{\left[#1\right]}

\newcommand{\field}[1]{\ensuremath{\mathbb{#1}}}
\newcommand{\abs}[1]{\left|#1\right|} 
\newcommand{\R}{\ensuremath{\field{R}}} 
\newcommand{\I}[1]{\ensuremath{\mathsf{1}_{\left\{#1\right\}}}} 
\newcommand{\PRP}[1]{\ensuremath{\mathsf{Pr}\left(#1\right)}} 
\newcommand{\ES}[1]{\ensuremath{\mathsf{E}\left[#1 \right]}} 

\renewcommand{\vec}[1]{\ensuremath{\boldsymbol{#1}}}

\newcommand{\CP}[2]{\ensuremath{\mathsf{Pr}\paren{\left.#1\right|#2}}}
\newcommand{\CD}[3]{\ensuremath{p_{#1}\paren{#2\left|#3\right.}}}
\newcommand{\CPMF}[3]{\ensuremath{p_{#1}\paren{#2\left|#3\right.}}}
\newcommand{\PMF}[2]{\ensuremath{p_{#1}\paren{#2}}}
\newcommand{\DSE}[1]{\ensuremath{\mathsf{H}\sqparen{#1}}}
\newcommand{\CDSE}[2]{\ensuremath{\mathsf{H}\sqparen{#1\left|#2\right.}}}
\newcommand{\MI}[2]{\ensuremath{\mathsf{I}\sqparen{#1;#2}}}
\newcommand{\KLD}[2]{\ensuremath{\mathsf{D}\sqparen{#1\left\|#2\right.}}}
\newcommand{\Null}[1]{\ensuremath{\mathsf{Null}\paren{#1}}}

\title{Privacy-aware Minimum Error Probability Estimation: An Entropy Constrained Approach}
\author{Ehsan Nekouei, Henrik  Sandberg, Mikael Skoglund and Karl H. Johansson 
\thanks{The authors ate with the  School of Electrical Engineering and Computer Science, KTH Royal Institute of Technology, Stockholm, Sweden. \{nekouei,hsan,skoglund, kallej\}@kth.se. This work is supported by the Knut and Alice Wallenberg Foundation, the
Swedish Foundation for Strategic Research, the Swedish Research Council.}}

\begin{document}
\maketitle
\thispagestyle{empty}

\begin{abstract}
This paper studies the design of an optimal privacy-aware estimator of a public random variable based on noisy measurements which contain private information. \textcolor{black}{The public random variable carries non-private information, however, its estimate will be correlated with the private information due to the estimation process. It is assumed that the  estimate of the public random variable is revealed to an untrusted party.} The objective is to design an estimator for the public random variable such that the leakage of the private information, via the estimation process, is kept below a certain level. The privacy metric is defined as the discrete conditional entropy of the private random variable, which carries the private information, given the output of the estimator. A binary loss function is considered for the estimation of the public random variable. It is shown that the optimal privacy-aware estimator is the solution of a (possibly infinite-dimensional) convex optimization problem. We next study the optimal perfect privacy estimator which ensures that the estimate of the public random variable is independent of the private information. \textcolor{black}{A necessary and sufficient condition is derived which guarantees that an estimator satisfies the perfect privacy requirement.} It is shown that the optimal perfect privacy estimator is the solution of a linear optimization problem. A sufficient condition for the existence of the optimal perfect privacy estimator is derived when the estimator has access to discretized measurements of sensors.
\end{abstract}

\section{Introduction}
\subsection{Motivation}
Networked systems play major roles in our society by providing critical services such as intelligent transportation and the smart grid. The operation of networked systems relies on the \emph{estimation process} wherein the objective is to obtain accurate values of certain variables based on noisy information collected by a set of sensors. However, the sensors' measurements not only contain information about the desired variable but also contain information which might be considered as private, \emph{e.g.,} information regarding stochastic events or unpredictable disturbances occurring in the sensors' environment. Hence, the output of an estimator of the desired variable may contain private information  due to the dependency of the sensors' measurements on the private variables. Thus, revealing the output of an estimator to an untrusted party might result in the loss of privacy. In what follows, we refer to the leakage of private information due to the estimation process as the \emph{privacy loss}.

 As a motivating example, consider a smart building application wherein sensors regularly collect noisy information about the temperature level in different areas of the building. The temperature measurement in each area of the building depends on the number of people, \emph{i.e,} the occupancy level, in that area which may be considered as private information. Thus, revealing an estimate of the temperature to an untrusted party, for storage or management purposes,  might result in the loss of privacy as the released data contain information regarding the occupancy level in different areas of the building. 

Due to the distributed structure of networked systems, the output of the estimation process is usually shared with untrusted parties, \emph{e.g.,} with a cloud-based storage system or with a cloud-based controller for temperature regulation in a smart building application. 
We refer to the untrusted party with access to the output of the estimator as the ``user.
Thus, to ensure the privacy of a networked system,  it is important to design privacy-aware estimators which provide accurate estimates of the desired variables based on the sensors' measurements and simultaneously  ensure that the leakage of private information due to the estimation process is kept below a certain level. 
\subsection{Contributions}
In this paper, we consider an estimation problem in which the sensors' measurements contain noisy information about a private random variable and a public random variable. It is assumed that the estimate of the public variable is revealed to the user. 
 Our objective is to design the optimal randomized estimator of the public random variable subject to a constraint on the privacy level of the private random variable.  The notion of conditional discrete entropy is used to quantify the leakage of the private information due to the estimation procedure. That is, the conditional discrete entropy of the private random variable given the output of the estimator  is considered as the privacy metric. The privacy metric captures the uncertainty of the user regarding the private random variable after observing the estimate of the public random variable. 

We first consider a single sensor estimation problem. It is shown that the optimal privacy-aware estimator is the solution of an infinite dimensional convex optimization problem when the estimator has access to the sensor's measurement, modeled as a continuous random variable. When the estimator receives a discretized version of the sensor's measurement, the optimal estimator design problem becomes a finite dimensional convex optimization problem. Necessary and sufficient optimality conditions for the optimal estimator are derived in the finite dimensional problem.

We also consider the optimal perfect privacy estimator design problem which ensures that the estimate of the public random variable is independent of the private random variable. When the estimator has access to the continuous measurements, we show that the  optimal perfect privacy estimator is the solution of an infinite dimensional linear optimization problem. It is shown that the feasible set of this optimization problem is always non-empty. When the estimator operates based on the discretized sensor's measurement, the optimal perfect privacy estimator is the solution of a finite dimensional linear optimization problem. It is shown that if the dimension of the null space of a certain matrix is non-zero, the feasible set of the optimal perfect privacy estimator design problem is non-empty. We also discuss the extension of the estimator design problem to a multi-sensor scenario.  

\subsection{Related Work}
The privacy aspect of hypothesis testing problems with a private and a public hypothesis has been studied in the literature, and various privacy-preserving solutions for improving the privacy level of hypothesis test problems have been proposed,  \emph{e.g.,} see \cite{HTS16, ST16, HT17, LSTC16}. In \cite{LO15}, the authors considered a hypothesis test problem with multiple sensors in which an eavesdropper intercepts the local decisions of a subset of sensors.  They studied the optimal decision rule  minimizing the Bayes risk at a fusion center subject to a privacy constraint at the eavesdropper.  In \cite{LO17}, the authors considered a similar set-up to that of \cite{LO15} and studied the optimal privacy-aware  Neyman-Pearson test with a private hypothesis.  The privacy of electricity consumers against an eavesdropper using demand management techniques and storage devices was studied in \cite{Li17}. 

   Privacy  preserving filters, for the state privacy problem in a cloud-based control application, were studied in \cite{TSSJ17} using the notion of directed information as the privacy metric. The privacy-aware controller design problem for a private Markov decision process problem in presence of an eavesdropper, with access to the input and output of the process, was studied in \cite{V13}. The authors in \cite{MO18} studied the privacy filter design problem for a public Markov chain, correlated with a private Markov chain, when both private and public chains are directly observable.
   
  \textcolor{black}{ The notion of differential privacy has been used to study the privacy-aware estimation, filtering and average consensus problems. The authors in \cite{NP14} proposed a filtering scheme for preserving the privacy of states or measurements of dynamical systems using the notion of differential privacy. The state estimation problem in a distribution power network subject to  differential privacy constraints for the consumers was studied in \cite{SDT15}. The authors in \cite{WHMD17} considered a  distributed multi-agent control problem and proposed a differential privacy scheme for preserving the privacy of the initial state as well as the preferred target way-points of each agent. Privacy-aware average consensus algorithms, for preserving the privacy of initial states of different agents, have been proposed in \cite{NTC17} \cite{MM17}. } 

 Information-theoretic methods for improving data privacy have been investigated in the literature, \emph{e.g.,} see \cite{KSK17}, \cite{BWI16}, \cite{CF12}, \cite{MS15} and references therein. In this line of research, the objective is to design privacy preserving filters which operate on a (directly observable) public random variable which is correlated with a private random variable.  The privacy filter is designed such that the distortion between the public variable and its processed version is minimized while a certain level of privacy is guaranteed. The current manuscript is different from this line of research in that, in our set up,  neither the public nor the private random variables are directly observable, and the sensors' observations contain noisy information about the public and private random variables. Moreover, in an estimation problem, one is interested in the true value of a variable based on a noisy observation rather than obtaining a low distortion representation of a directly observable random variable.

 \textcolor{black}{In \cite{AAL16}, the authors considered the problem of adding stochastic distortion to a public variable, which depends on a private information, such that $(i)$ the mean square error (MSE) of recovering the original variable from its  distorted version is minimized, $(ii)$ the minimum MSE of recovering the private information from the distorted variable stays above a certain level. Their results were extend in \cite{ADAL17} under the Hamming distance as the distortion criterion and the efficiency of these methods was analysed in \cite{ADAL17j}.}

Perfect privacy filters in the context of data privacy have been studied in \cite{RG17}. The results of \cite{RG17} are mainly derived based on the assumption that the public random variable is directly observable and takes finite values. Also, the perfect privacy condition in \cite{RG17} requires characterization of the extreme points of a certain convex polytope which depends on the null space of a probability transition matrix. Finally, we note that the relation between the perfect privacy condition  and the notion of maximal correlation has been studied in \cite{CMMV17}.

 Different from the existing work, we consider the problem of recovering a public random variable from a noisy observation which depends on both private and public random variables. That is, in our set-up, the sensor does not have (direct) access to the private or public random variables. Our objective is to design an estimator for the public random variable with a guarantee on the leakage of the private information. Our work differs from \cite{RG17}, in that,  we assume that the sensor's observations take values in $\R$ which results in an infinite-dimensional convex optimization problem. In our work, the optimal perfect privacy estimator  with discretized measurements is obtained by solving a linear optimization problem which does not require finding the extreme points of a convex polytope. Moreover, we study the  perfect privacy condition for both continuous and discretized measurements. Our results show that the perfect privacy with the discretized measurements  depends on the null space of a matrix with positive and negative entries which is different from a transition probability matrix. 
\begin{figure*}
\centering{\includegraphics[scale=0.15]{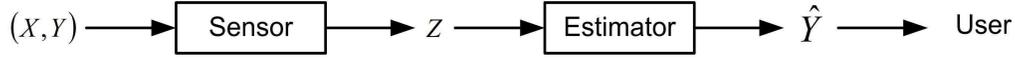}}
\caption{A single sensor estimation set up.}\label{Fig: F1}
\hrule
\end{figure*}

The rest of this paper is organized as follows. Our system model and assumptions are described in the next section. Section \ref{Sec: PSS} presents our results on the optimal privacy-aware and optimal perfect privacy estimator design problems with a single sensor. These results are extended to a multi-sensor set-up in Section \ref{Sec: PMS}.  Our numerical results are preseneted in Section \ref{Sec: NM}, followed by the concluding remarks in Section \ref{Sec: Conc}.  
\section{ System Model}\label{Sec: SM}
Consider an estimation problem with one sensor in which the sensors' measurement, denoted by $Z$, contains noisy information about two, possibly correlated,  discrete random variables $X$ and $Y$. The random variable $Y$ contains public information, and upon observing $Z$, an estimate of $Y$ is delivered to an untrusted party, hereafter, named the user. The random variable $X$ carries information which should be kept hidden from any untrusted party. The support sets of $X$, $Y$ and $Z$ are denoted by $\mathcal{X}=\left\{x_1,\cdots,x_n\right\}$, $\mathcal{Y}=\left\{y_1,\cdots,y_m\right\}$ and $\mathcal{Z}=\R$, respectively. Through this paper, we assume that the random variable $Z$ is absolutely continuous with respect to Lebesgue measure on $\R$ with the probability density function $p_Z\paren{z}$.  A pictorial representation of our system model is illustrated in Fig. \ref{Fig: F1}.

 Let $\hat{Y}\paren{Z}$ denote an estimate of $Y$ based on the observation $Z$. The user can infer information about the private random variable $X$ by observing $\hat{Y}\paren{Z}$ since $\hat{Y}\paren{Z}$ is usually correlated with $X$. Thus,  revealing $\hat{Y}\paren{Z}$ to the user may result in the privacy loss, \emph{i.e.,} the user can infer about $X$ by observing $\hat{Y}\paren{Z}$. In what follows, by privacy loss we mean the leakage of the private information due to the estimation task.

\begin{figure*}
\begin{eqnarray}\label{Eq: CDE}
\CDSE{X}{\hat{Y}_P\paren{Z}}=-\sum_{y\in\mathcal{Y}}\PRP{\hat{Y}_P\paren{Z}=y}\sum_{x\in\mathcal{X}}\CP{X=x}{\hat{Y}_P\paren{Z}=y}\log\CP{X=x}{\hat{Y}_P\paren{Z}=y}
\end{eqnarray}
\hrule
\end{figure*}
\begin{figure*}
\begin{align}\label{Eq: PC}
\CDSE{X}{\hat{Y}_P\paren{Z}}&=\DSE{X}-\sum_j\PRP{X=x_j}\KLD{\CPMF{\hat{Y}_P}{y}{X=x_j}}{\PMF{\hat{Y}_P}{y}}
\end{align}
\hrule
\end{figure*}

In this paper, our objective is to design an estimator for the public random variable $Y$ which minimizes a desired loss function while the information leakage about the private variable $X$ is kept below a certain level. An estimator of $Y$ is a (possibly randomized) map from $\mathcal{Z}$ to $\mathcal{Y}$. Let $P\paren{z}=\left[P_i\paren{z}\right]_{i=1}^m$ denote a set of positive functions  where $m=\abs{\mathcal{Y}}$ and $P_i\paren{z}$  is defined on the support set of $\mathcal{Z}$ with $\sum_{i=1}^mP_i\paren{z}=1$ for all $z\in \mathcal{Z}$. Then, a randomized estimator of $Y$ can be expressed as 
\begin{eqnarray}\label{Eq: RE}
\hat{Y}_P\paren{z}=
\left\{
\begin{array}{cc}
 y_1   &\text{w.p.} \quad P_1\paren{z} \\  
 \vdots&\vdots\nonumber\\
  y_m   & \text{w.p.}\quad  P_m\paren{z}\\
\end{array}
 \right.\\
\end{eqnarray}
where w.p. stands for with probability. According to \eqref{Eq: RE}, if the sensor's measurement is equal to $z$, the estimator declares $y_i$ as the estimate of $Y$ with probability  $P_i\paren{z} $.

 We use a binary loss function, denoted by $L\paren{Y,\hat{Y}_P\paren{Z}}$, to quantify the estimation loss of the private random variable which is expressed as
\begin{align}
L\paren{Y,\hat{Y}_P\paren{Z}}=
\left\{
\begin{array}{cc}
 1   & Y\neq \hat{Y}_P\paren{Z} \\  
 0   & Y= \hat{Y}_P\paren{Z}\nonumber 
\end{array}
 \right.
\end{align}
Thus, the estimation loss is equal to 1 if the output of the estimator is different from the true value of $Y$ and there is no loss if these two values agree.\textcolor{black}{
\subsection{Motivating Examples}
In this subsection, we provide three motivating examples of privacy-aware estimation problem which include public and private random variables.
\subsubsection{Temperature and ${\rm CO}_2$ Estimation In A Building}
Consider a smart building application in which the objective is to estimate the temperature and ${\rm CO}_2$ levels inside a building based on the noisy measurements of sensors. Here, the random variable $X$ represents the number of occupants in the building which is considered as private information, $Y$ represents the actual temperature and ${\rm CO}_2$ levels and $Z$ represents noisy measurements of temperature and ${\rm CO}_2$ obtained from sensors. In this application, $Y$ depends on $X$ directly. The objective is to provide an accurate estimate $Y$, based on $Z$, for an untrusted user, \emph{e.g.,} a cloud-based controller, while the dependency of $\hat{Y}$ on $X$ is kept below a certain level.
\subsubsection{Smart Meter Problem}
Consider a smart grid application wherein the objective is to estimate the number of active smart appliances in a household based on its electricity consumption level. In this application, the absence/presence of tenants can be modeled using the binary random variable $X$, which contains private information, $Y$ represents the number of active smart appliances and  $Z$ represents the electricity consumption level. In this application, an estimate of $\hat{Y}$ is provided to an untrusted party, \emph{e.g.,} a utility company. The objective is to design an estimator of $Y$ such that the leakage of the private information, after $\hat{Y}$ is revealed to the untrusted party, is kept below a certain level.    
\subsubsection{Counting The Vehicles On A Road} Consider the problem of counting vehicles on a road based on the position estimates of individual vehicles. In this application, $X$ represents the velocity of a vehicle which is considered as its private information, $Y$ represents the actual position of the vehicle and $Z$ denotes the noisy measurement of the position of the vehicle. The estimates of position by individual vehicles can used by a traffic operator, \emph{i.e.,} an untrusted party, to estimate the number of vehicles on the road at a given time.  In this application, the objective is to design an estimator of $Y$ based on $Z$ for individual vehicles such that leakage information about each vehicle's velocity, after revealing its corresponding position estimate $\hat{Y}$ to the traffic operator, is kept below a certain level.
}
\subsection{Privacy Metric}
In this paper, we consider the conditional discrete entropy, or equivocation,  as the privacy metric. The conditional discrete entropy of $X$ given $\hat{Y}_P\paren{Z}$, denoted by $\CDSE{X}{\hat{Y}_P\paren{Z}}$, is defined in \eqref{Eq: CDE}. Our choice of privacy metric is motivated by the fact that $\CDSE{X}{\hat{Y}_P\paren{Z}}$ captures the ambiguity of the user about $X$ after observing $\hat{Y}_P\paren{Z}$. Thus, the privacy loss decreases as $\CDSE{X}{\hat{Y}_P\paren{Z}}$ becomes large since the user becomes more uncertain about the value of $X$ as $\CDSE{X}{\hat{Y}_P\paren{Z}}$ increases. Since conditioning reduces entropy \cite{CT06}, we have 
\begin{align}
0\leq \CDSE{X}{\hat{Y}_P\paren{Z}}\leq \DSE{X}\nonumber
\end{align}
which implies that the maximum privacy is achieved if $\CDSE{X}{\hat{Y}_P\paren{Z}}=\DSE{X}$. Recall that if $X$ and $\hat{Y}_P\paren{Z}$ are independent, $\hat{Y}_P\paren{Z}$ contains no information about $X$ and the user has maximum ambiguity about $X$ after observing $\hat{Y}_P\paren{Z}$, \emph{i.e.,} $\CDSE{X}{\hat{Y}_P\paren{Z}}=\DSE{X}$.

The other motivation for the choice of the privacy metric in this paper is the fact that the error probability of estimating $X$ after observing $\hat{Y}_P\paren{Z}$ can be lower bounded in terms of $\CDSE{X}{\hat{Y}_P\paren{Z}}$ using Fano's inequality \cite{CT06}:
\begin{align}
 \PRP{X\neq\hat{X}\paren{\hat{Y}}}\geq \frac{\CDSE{X}{\hat{Y}_P\paren{Z}}-1}{\log\abs{\mathcal{X}}}
 \end{align}
 where $\hat{X}\paren{\hat{Y}}$ is an arbitrary estimator of $X$ and  $\abs{\mathcal{X}}$ is the cardinality of the support set of $X$. Note that this lower bound is independent of the estimator of $X$, \emph{i.e.,} it holds for all possible estimators. Thus, by adjusting the value of $\CDSE{X}{\hat{Y}_P\paren{Z}}$, a desired privacy level of the private random variable can be guaranteed as long as $\abs{\mathcal{X}}>2$. 

\section{Privacy-aware Optimal Estimation Problem: Single Sensor Case}\label{Sec: PSS}
In this section, the design of the optimal privacy-aware estimator of the public random variable is studied. We first study the privacy-aware and the perfect privacy estimator design problems when the estimator has access to the sensor's measurement. We then investigate the estimator design problems when the discretized sensor's measurements are available at the estimator.
 \subsection{Optimal Estimator With Continuous Measurements}
 In this subsection, we assume that the estimator has access to $Z$ which is a continuous random variable. In this case, the optimal design of the estimator subject to the privacy constraint is given by the solution of the following optimization problem:
\begin{align}\label{Eq: OP}
\underset{ \left\{P_i\paren{z}\right\}_{i=1}^m }{\rm minimize} \quad& \ES{L\paren{Y,\hat{Y}_P\paren{Z}}}\nonumber\\
& P_i\paren{z}\geq 0, \forall i\nonumber\\
&\sum_iP_i\paren{z}=1, \quad \forall z\nonumber\nonumber\\
& \CDSE{X}{\hat{Y}_P\paren{Z}}\geq \mathsf{H}_0
\end{align}

According to this optimization problem, the functions $\left[P_i\paren{z}\right]_i$ are chosen such that the average estimation loss is minimized and, a certain level of privacy is ensured by keeping the conditional discrete entropy of $X$ given $\hat{Y}_P\paren{Z}$ above the desired level $\mathsf{H}_0$. 

The optimization problem \eqref{Eq: OP} is a functional optimization problem defined on the space of bounded measurable functions from $\R$ to $\R$, \emph{i.e.,} $\mathsf{B}\paren{\R,\R}$. Note that $\mathsf{B}\paren{\R,\R}$ forms a Banach space under the supremum norm and $P_i\paren{z}$ belongs to the cone of positive functions in $\mathsf{B}\paren{\R,\R}$ for all $i$. Next lemma derives an expression for  the objective function of the optimization problem \eqref{Eq: OP}. 
\begin{lemma}\label{Lem: Obj}
The objective function in \eqref{Eq: OP} can be written as 
\begin{align}
 1-\sum_i\int P_i\paren{z}\CP{Y=y_i}{Z=z}p_Z\paren{z}dz\nonumber
\end{align}
where $p_Z\paren{z}$ is the probability density function of $Z$ and $\CP{Y=y_i}{Z=z}$ is the conditional probability of the event $Y=y_i$ given the event $Z=z$.
\end{lemma}
\begin{IEEEproof}
See Appendix \ref{App: Obj}.
\end{IEEEproof}
According to Lemma \ref{Lem: Obj}, the objective function of the optimization problem \eqref{Eq: OP} is  linear in the decision variables $P\paren{z}=\left[P_i\paren{z}\right]_i$. 
Next lemma studies the convexity of the privacy constraint. 
\begin{lemma}\label{Lem: P-Conv}
The privacy constraint can be written as \eqref{Eq: PC} where $\DSE{X}$ is the discrete entropy of $X$, $\PMF{\hat{Y}_P}{y}$ and  $\CPMF{\hat{Y}_P}{y}{X=x_j}$ denote the probability mass function of $\hat{Y}_P\paren{Z}$ and the conditional probability mass function of $\hat{Y}_P\paren{Z}$ given $X=x_j$, respectively, and $\KLD{\cdot}{\cdot}$ denotes Kullback-Leibler (KL) divergence (relative entropy). Furthermore, the privacy constraint is convex in $P\paren{z}$.
\end{lemma}

\begin{IEEEproof}
See Appendix \ref{App: P-Conv}.
\end{IEEEproof}
The objective function in the optimization problem \eqref{Eq: OP} is linear and the constraint set is convex. Thus, the optimization problem \eqref{Eq: OP} is a convex optimization problem. This result is formally stated in the next theorem.
\begin{theorem}
The optimal privacy-aware estimator of the public random variable can be designed by solving the convex optimization problem  \eqref{Eq: OP}. 
\end{theorem}

\subsection{Optimal Estimator With Discretized  Measurements}
In this subsection, we assume that the estimator has only access to a discretized version of the sensor's measurement. This allows us to express the privacy-aware estimator design problem as a finite dimensional convex optimization problem. To this end, let $\left\{B_i\right\}_{i=1}^N$ denote a partition of $\R$ where $B_1$ and $B_N$ are semi-infinite intervals and where $B_i$, $2\leq i\leq N-1$, are of the form $B_i=\left[a_{i-1},a_i\right]$, $a_i>a_{i-1}$. A randomized estimator of $Y$, based on the discretized measurements, is defined as 
\begin{eqnarray}
\hat{Y}_P\paren{z}=
\left\{
\begin{array}{ccc}
 y_1   &\text{w.p.} \quad P_{1l}, & \text{if} \quad z\in B_l\nonumber\\  
 \vdots&\vdots&\vdots\nonumber\\
  y_m   & \text{w.p.}\quad  P_{ml},& \text{if} \quad z\in B_l\nonumber\\
\end{array}
 \right.\nonumber
\end{eqnarray}
where $\sum_{i}P_{il}=1$ for all $l\in\left\{1,\cdots,N\right\}$. Thus, the estimator selects $y_i$ as its output with probability $P_{il}$ if the sensor's measurement belongs to $B_l$. The optimal privacy-aware estimator, under discretized measurements, is given by 
\begin{align}\label{Eq: OP-D}
\underset{ \left\{P_{il}\right\}_{i,l} }{\rm minimize} \quad& \ES{L\paren{Y,\hat{Y}_P\paren{Z}}}\nonumber\\
& P_{il}\geq 0, \forall i, l\nonumber\\
&\sum_iP_{il}=1, \quad \forall l\nonumber\nonumber\\
& \CDSE{X}{\hat{Y}_P\paren{Z}}\geq \mathsf{H}_0
\end{align}

In this case, the probability of correct estimation can be expressed as
\begin{align}
&\PRP{Y=\hat{Y}_P\paren{Z}}\nonumber\\
&=\sum_l\PRP{Y=\hat{Y}_P\paren{Z},Z\in B_l}\nonumber\\
&=\sum_{il}\CP{Y=y_i}{Z\in B_l}\CP{Z\in B_l}{Y=y_i}\PRP{Y=y_i}\nonumber\\
&=\sum_{il}P_{il}\CP{Z\in B_l}{Y=y_i}\PRP{Y=y_i}\nonumber
\end{align}

Moreover, $\PMF{\hat{Y}_P}{y}$ and  $\CPMF{\hat{Y}_P}{y}{X=x_j}$ can be written as 
\begin{align}
\PMF{\hat{Y}_P}{y_i}&=\PRP{\hat{Y}_P\paren{Z}=y_i}\nonumber\\
&=\sum_l\CP{\hat{Y}_P\paren{Z}=y_i}{Z\in B_l}\PRP{Z\in B_l}\nonumber\\
&=\sum_{l} P_{il}\PRP{Z\in B_l}
\end{align}
and 
\begin{align}
&\CPMF{\hat{Y}_P}{y_i}{X=x_j}\nonumber\\
&=\CP{\hat{Y}_P\paren{Z}=y_i}{X=x_j}\nonumber\\
&=\sum_l\CP{\hat{Y}_P\paren{Z}=y_i}{Z\in B_l}\CP{Z\in B_l}{X=x_j}\nonumber\\
&=\sum_lP_{il}\CP{Z\in B_l}{X=x_j}\nonumber
\end{align}
Following the proofs of Lemmas \ref{Lem: Obj} and \ref{Lem: P-Conv}, it can be shown that the optimization problem \eqref{Eq: OP-D} is a convex optimization problem. 

The next lemma states the KKT necessary and sufficient optimality conditions for the optimization problem \eqref{Eq: OP-D}.

 \begin{lemma}\label{Lem: KKT}
Let $\vec{P}^\star$ denote the optimal solution of the optimization problem \eqref{Eq: OP-D}. Then, we have \eqref{Eq: KKT-Opt}  where $\mu^\star$ is the dual optimal variable associated with the privacy constraint and $\vec{\lambda}^\star$ is the vector of  dual optimal variables associated with the equality constraints.
  \end{lemma}
  
\begin{figure*}
\begin{align}\label{Eq: KKT-Opt}
&\mu^\star\log\prod_j\paren{\frac{\sum_l P^\star_{hl}\CP{Z\in B_l}{X=x_j}}{\sum_{l}P^\star_{hl}\PRP{Z\in B_l}}}^{\PRP{Z\in B_k,X=x_j}}\left\{
\begin{array}{cc}
={\CP{Z\in B_k}{Y=y_h}\PRP{Y=y_h}+\lambda_k^\star} \quad \text{if}&  0<P^\star_{hk}<1\\
\leq  {\CP{Z\in B_k}{Y=y_h}\PRP{Y=y_h}+\lambda_k^\star}  \quad\text{if}& P^\star_{hk}=0\\
\geq {\CP{Z\in B_k}{Y=y_h}\PRP{Y=y_h}+\lambda_k^\star} \quad \text{if}& P^\star_{hk}=1
\end{array}
\right.\nonumber\\
&\sum_hP_{hl}^\star=1 \quad \forall l\nonumber\\
&\mu^\star\geq 0\nonumber\\
& \CDSE{X}{\hat{Y}_{\vec{P}^\star}\paren{Z}}\geq \mathsf{H}_0 \nonumber\\
&\mu^\star\paren{\mathsf{H}_0 -\CDSE{X}{\hat{Y}_{\vec{P}^\star}\paren{Z}}}=0
\end{align}
\hrule
\end{figure*}

\begin{IEEEproof}
See Appendix \ref{App: KKT}.
\end{IEEEproof}

\subsection{Optimal Perfect Privacy Estimator With Discretized Measurements}\label{Subsec: PP-DM}
In this subsection, we first define the perfect privacy condition, and derive a necessary and sufficient condition for an estimator to satisfy the perfect privacy requirement. Then, we show that the optimal perfect privacy estimator can be obtained by solving a linear optimization problem.
 \begin{definition}
 An estimator of the public random variable satisfies the perfect privacy condition if the output of the estimator is independent of the private random variable.
 \end{definition}

 Before proceeding with the derivation of the perfect privacy condition, we first define the matrix $\Phi$ as 
\begin{eqnarray}\label{Eq: Phi}
\Phi=\left[
\begin{array}{ccc}
\phi_{11}&\cdots&\phi_{1N}\nonumber\\
\vdots& &\vdots\nonumber\\
\phi_{n1}&\cdots&\phi_{nN}\nonumber\\
\end{array}
\right]
\end{eqnarray}
where $\phi_{jl}=\CP{Z\in B_l}{X=x_j}-\PRP{Z\in B_l}$. Also, the vector $\vec{P}_i$ is defined as 
\begin{eqnarray}
\vec{P}_i=
\left[
\begin{array}{c}
P_{i1}\nonumber\\
\vdots\nonumber\\
P_{iN}\nonumber\\
\end{array}
\right]
\end{eqnarray}
which is the collection of randomization probabilities associated with selecting $y_i$, as the output of an estimator, for different bins. Next lemma derives a necessary and sufficient condition for the randomization probabilities which ensure the prefect privacy for the estimation.
\begin{lemma}\label{Lem: PP-C-CM}
An estimator satisfies the perfect privacy condition if and only if $\vec{P}_i\in \Null{\Phi}$ for all $i\in\left\{1,\cdots,m\right\}$ where $\Null{\Phi}$ is the null space of the matrix $\Phi$ defined in \eqref{Eq: Phi}. 
\end{lemma}

\begin{IEEEproof}
See Appendix \ref{App: PP-C-DM}.
\end{IEEEproof}

According to this lemma, a randomized estimator satisfies the perfect privacy condition if the vector of randomization probabilities associated with each element of $\mathcal{Y}$, \emph{i.e.,} the support set of the public random variable, lies in the null space of the matrix $\Phi$. \textcolor{black}{We note that the perfect privacy conditions have been studied in the literature for different settings, \emph{e.g.,} see \cite{RG17} and \cite{AAL14}. For example, the perfect privacy conditions in \cite{RG17} require that the dimension of the null space of certain transition probability matrices to be non-zero. However, the perfect privacy condition in our paper requires that the randomization probability vectors $\vec{P}_i$s to lie in the null space of the matrix $\Phi$ which is not a transition probability matrix.}

We next study the optimal perfect privacy estimator design problem, \emph{i.e.,} finding the optimal estimator among all the estimators which satisfy the perfect privacy condition. The optimal perfect privacy estimator is given by the solution of 
\begin{align}\label{Eq: OP-PP}
\underset{ \left\{P_{il}\right\}_{i,l} }{\rm minimize} \quad& \ES{L\paren{Y,\hat{Y}_P\paren{Z}}}\nonumber\\
& P_{il}\geq 0, \forall i, l\nonumber\\
&\sum_iP_{il}=1, \quad \forall l\nonumber\nonumber\\
& \Phi\vec{P}_i=\vec{0}, \forall i
\end{align}

Note that the perfect privacy condition in the optimization problem above is linear in the decision variables, \emph{i.e.,} the randomization probabilities. Thus, the optimal perfect privacy estimator can be obtained by solving a linear optimization problem. 

Let $\mathcal{P}=\mathcal{P}_1\times\cdots\times\mathcal{P}_N$ denote the joint probability simplex corresponding to the randomization probabilities of bins, \emph{i.e.,} for each bin $l$, $l\in\left\{1,\cdots,N\right\}$, we have $\paren{P_{1l},\cdots,P_{ml}}^\top\in\mathcal{P}_l$. The feasible set of the optimization problem above is the intersection of the set $\mathcal{P}$ and the perfect privacy condition. We next derive a sufficient condition which ensures that the feasible set of the optimization problem \eqref{Eq: OP-PP}  is non-empty.
\begin{lemma}\label{Lem: PP-Feas}
If the number of discretization bins, \emph{i.e.,} $N$, is more than the size of the support set of the private random variable $X$, then the feasible set of \eqref{Eq: OP-PP} is non-empty.
\end{lemma}

\begin{IEEEproof}
See Appendix \ref{App: PP-DM-Feas}.
\end{IEEEproof}

\subsection{Optimal Perfect Privacy Estimator with Continuous Measurements}
\textcolor{black}{In this subsection, the perfect privacy condition with continuous measurements is derived and it is  shown that the optimal perfect privacy estimator is the solution of an infinite dimensional linear optimization problem.}

According to the perfect privacy definition, the output of the estimator, $\hat{Y}_P\paren{Z}$, satisfies the perfect privacy requirement if and only if we have 
\begin{align}
\PRP{\hat{Y}_P\paren{Z}=y_i,X=x_j}=\PRP{\hat{Y}_P\paren{Z}=y_i}\PRP{X=x_j}\nonumber
\end{align} 
for all $i,j$. Note that, $\PRP{\hat{Y}_P\paren{Z}=y_i,X=x_j}$ and $\PRP{\hat{Y}_P\paren{Z}=y_i}$ can be written as 
\begin{align}
&\PRP{\hat{Y}_P\paren{Z}=y_i,X=x_j}\nonumber\\
&\hspace{1cm}=\int\CP{\hat{Y}_P\paren{Z}=y_i,X=x_j}{Z=z}p_Z(z)dz\nonumber\\
&\hspace{1cm}=\PRP{X=x_j}\int P_i\paren{z}\CD{Z}{z}{X=x_j}dz\nonumber
\end{align}
and 
\begin{align}
\hspace{-2cm}\PRP{\hat{Y}_P\paren{Z}=y_i}=\int P_i\paren{z}p_Z\paren{z}dz\nonumber
\end{align}
Thus, for $\PRP{X=x_j}\neq 0$, the perfect privacy condition can be written as 
\begin{align}
\int P_i\paren{z}\paren{p_Z\paren{z}-\CD{Z}{z}{X=x_j}}dz=0\nonumber
\end{align}

Thus, the optimal perfect privacy estimator design problem with continuous measurements can be written as 
 \begin{align}\label{Eq: OP-PP-CM}
 \underset{ \left\{P_i\paren{z}\right\}_{i=1}^m }{\rm minimize} \quad& \ES{L\paren{Y,\hat{Y}_P\paren{Z}}}\nonumber\\
 &\hspace{-1.5cm} P_i\paren{z}\geq 0, \forall i\nonumber\\
 &\hspace{-1.5cm}\sum_iP_i\paren{z}=1, \quad \forall z\nonumber\nonumber\\
 &\hspace{-1.5cm}\int P_i\paren{z}\paren{p_Z\paren{z}-\CD{Z}{z}{X=x_j}}dz=0 \quad \forall i,j 
 \end{align}
 Note that perfect privacy constraint in the optimization problem above is a linear constraint.
 
We next show that the feasible set of the optimal perfect privacy estimation problem, in the continuous case, is always non-empty.
\begin{lemma}\label{Lem: PP-CM-Feas}
The feasible set of the optimal perfect privacy estimator design problem with continuous measurements in \eqref{Eq: OP-PP-CM} is always non-empty.
\end{lemma}
\begin{IEEEproof}
See Appendix \ref{App: PP-CM-Feas}.
\end{IEEEproof}
\section{Extension to Multi-sensor Case With Discrete Measurements}\label{Sec: PMS}
In this section, we extend the privacy-aware and perfect privacy estimation problems to a multi-sensor scenario. To this end, consider a multi-sensor estimation problem with $M$ sensors. Let $Z_s$, $s\in\left\{1,\cdots, M\right\}$, denote the measurement of sensor $s$ which depends on the private random variable $X$ and the public random variable $Y$. We assume that the sequence of random variables $\left\{Z_s\right\}_s$ are conditionally independent given $X$ and $Y$, \emph{i.e.,} the joint distribution of  $\left\{Z_s\right\}_s$ given $X$ and $Y$ is equal to the product of the conditional distributions of the sensors' measurements given $X$ and $Y$.  We only focus on the optimal estimator design with discretized measurements. These results can be easily extended to the continuous measurements case.
\begin{figure*}
\centering{\includegraphics[scale=0.15]{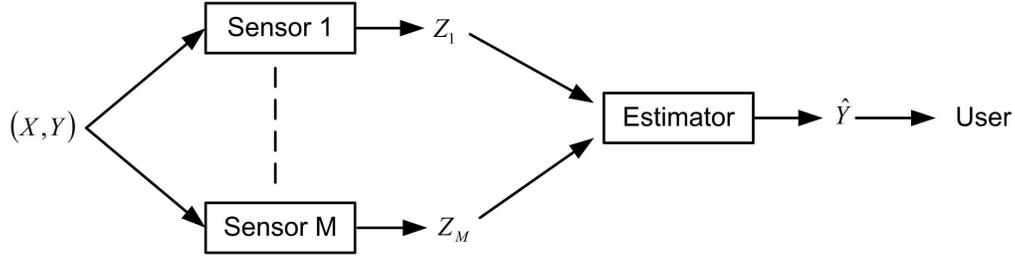}}
\caption{The multi-sensor estimation set up.}\label{Fig: FMS}
\end{figure*}

In the multi-sensor case, each sensor $s$ transmits the index of the bin, corresponding to the realization of its measurement $Z_s$, to a (common) remote estimator. Fig. \ref{Fig: FMS} illustrates the multi-sensor estimation set-up.  Let $\tilde{Z}_s$ denote the information transmitted to the estimator by sensor $s$. Then, we have $\tilde{Z}_s=l$ if $Z_s\in {\rm B}_l$ and the random variable $\tilde{Z}_s$ takes values in $\left\{1,\cdots,N\right\}$ for all $s$ where $N$ is the number of bins. The estimator computes an estimate of $Y$, based on the received information from sensors,  subject to the privacy constraint. More precisely, the  estimator first forms the vector $C=\left(C_1,\cdots,C_N\right)$ according to $C_l=\sum_i\I{\tilde{Z}_i=l}$ for $l\in\left\{1,\cdots,N\right\}$ where $C_l$ denotes the number of sensors which transmit $l$ as their bin indicies to the remote estimator. Note that the random vector $C$ takes values in $\mathcal{C}=\left\{\paren{c_1,\cdots,c_N}\in \mathbb{N}_0^N, \sum_l c_l=M\right\}$ where $\mathbb{N}_0$ is the set of non-negative integers. In what follows, we use $c=\paren{c_1,\cdots,c_N}\in \mathcal{C}$ to represent a realization of $C$.  

Given $C=c$, the output of a randomized estimator of $Y$ can be expressed as
\begin{eqnarray}
\hat{Y}_P\paren{c}=
\left\{
\begin{array}{ccc}
 y_1   &\text{w.p.} \quad P_{1c}, \\  
 \vdots&\vdots\nonumber\\
  y_m   & \text{w.p.}\quad  P_{mc},\\
\end{array}
 \right.\nonumber
\end{eqnarray}
where $P_{ic}\geq 0$ and $\sum_iP_{ic}=1$. That is, the estimator selects $y_i$ as the estimate of $Y$ with probability $P_{ic}$. Under  the binary loss function defined in Section \ref{Sec: SM}, the objective function of the estimation problem can be written as 
\begin{align}
1-\PRP{Y=\hat{Y}\paren{C}}=1-\sum_{i,c}P_{ic}\CP{C=c}{Y=y_i}\nonumber
\end{align}
which is linear in the decision variables $\left\{P_{ic}\right\}_{ic}$. Thus, the optimal privacy-aware estimator design problem can be written as 
\begin{align}
\underset{ \left\{P_{ic}\right\}_{i,c} }{\rm minimize} \quad & 1-\sum_{i,c}P_{ic}\CP{C=c}{Y=y_i}\nonumber \\
& P_{ic}\geq 0, \forall i, c\in\mathcal{C}\nonumber\\
&\sum_iP_{ic}=1, \quad \forall c\in\mathcal{C}\nonumber\nonumber\\
& \CDSE{X}{\hat{Y}_P\paren{C}}\geq \mathsf{H}_0\nonumber
\end{align}
Following the proof of Lemma \ref{Lem: P-Conv}, it is straightforward to show that the privacy constraint in the optimization problem above is a convex constraint. Thus, the optimal privacy-aware estimator design problem is a convex optimization problem. Moreover, the optimal perfect privacy estimator in the multi-sensor case be obtain by solving the following optimization problem
\begin{align}
\underset{ \left\{P_{ic}\right\}_{i,c} }{\rm minimize} \quad & 1-\sum_{i,c}P_{ic}\CP{C=c}{Y=y_i}\nonumber \\
&\hspace{-1cm} P_{ic}\geq 0, \forall i, c\in\mathcal{C}\nonumber\\
&\hspace{-1cm}\sum_iP_{ic}=1, \quad \forall c\in\mathcal{C}\nonumber\nonumber\\
&\hspace{-1cm}\sum_{c}P_{ic}\paren{\CP{C=c}{X=x_j}-\PRP{C=c}}=0,\forall j, i\nonumber
\end{align}
Following similar steps as those in the proof of Lemma \ref{Lem: PP-Feas}, one can show that the feasible set of the above  optimization problem is non-empty  if $\abs{\mathcal{C}}>n$. 

Finally, we derive an expression for the objective function of the estimator design problem when the sensors' measurements, \emph{i.e.,} $Z_s$s, are identically distributed. Note that we have $\CP{C=c}{Y=y_i}=\sum_{j} \CP{C=c}{Y=y_i,X=x_j}$. Moreover, using the fact that the sensors' measurements are conditionally independent given $X$ and $Y$,  we have  
\begin{align}
&\CP{C=c}{Y=y_i,X=x_j}\nonumber\\
&\hspace{2cm}=\frac{M!}{c_1!\cdots c_N!} \prod_l\CP{Z_s\in B_l}{Y=y_i,X=x_j}^{c_l}\nonumber
\end{align} 
for identically distributed measurements.
\section{Numerical Results}\label{Sec: NM}
In this section, we consider a multi-sensor estimation problem with $M$ sensors. It is assumed that the measurement of sensor $s$ is given by $Z_s=0.6Y+0.4X+N_s$ where $X,Y\in\left\{0,1\right\}$, $\PRP{X=0}=0.7$, $\PRP{Y=0}=0.5$, and $N_s$ is a zero mean Gaussian random variable with variance equal to $0.01$. The random variable $0.4X+N_s$ represents the measurement noise of sensor $s$. Thus, when the value of $X$ changes from $0$ to $1$, the mean of measurement noise changes from $0$ to $0.4$. We assume that the change in the mean of the measurement noise is triggered by an external stimulant, \emph{e.g.,} a change in the temperature, affects all the sensors simultaneously and it should be kept private. The sequence of random variables $\left\{N_s\right\}_s$ is assumed to be independent and identically distributed. In our numerical results, the number of discretization bins is equal to $4$ for all sensors. 

\begin{figure}[!htp]
\centering{\includegraphics[scale=0.6]{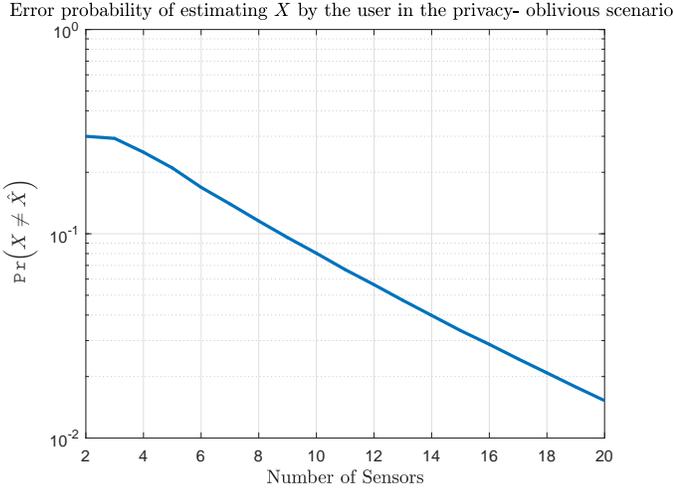}}
\caption{The error probability of the optimal estimator of $X$ in the privacy-oblivious scenario as a function of the number of sensors.}\label{Fig: F2}
\end{figure}

We compare the performance of the optimal perfect privacy estimator with that of a privacy-oblivious estimation scenario. In the privacy-oblivious scenario, each sensor individually estimates $Y$ (with no privacy constraint) and the user receives the estimates of $Y$ from sensors. Let $\hat{Y}_s^{\rm po}$ denote the optimal estimate of $Y$ by sensor $s$ based on its local measurement $Z_s$ with no privacy constraint.  Then, under the privacy oblivious scenario, the user can estimate both $X$ and $Y$ using the received local estimates of sensors, \emph{i.e.,} $\left\{\hat{Y}_s^{\rm po}\right\}_s$. Fig. \ref{Fig: F2} shows the error probability of estimating $X$ using $\left\{\hat{Y}_s^{\rm po}\right\}_s$ as a function of the number of sensors. As this figure shows, the ambiguity of the user regarding the private variable $X$ decreases as the number of sensors becomes large. This observation suggests that revealing the local estimates of sensors to an untrusted user results in a high level of leakage of private information. 

\begin{figure}[!htp]
\subfigure[]
{\centering{\includegraphics[scale=0.6]{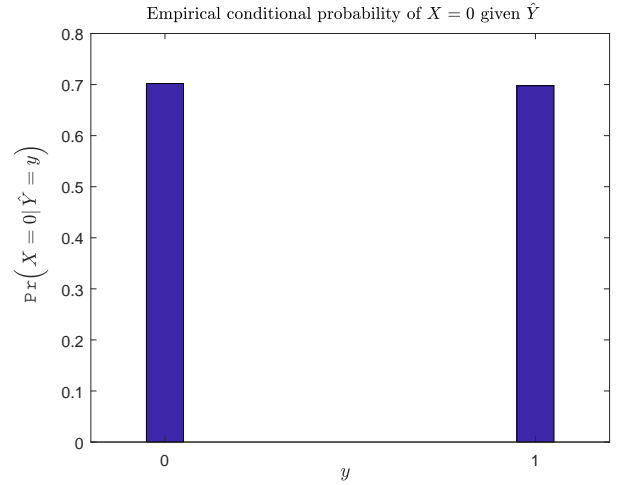}}}
\subfigure[]
{\centering{\includegraphics[scale=0.6]{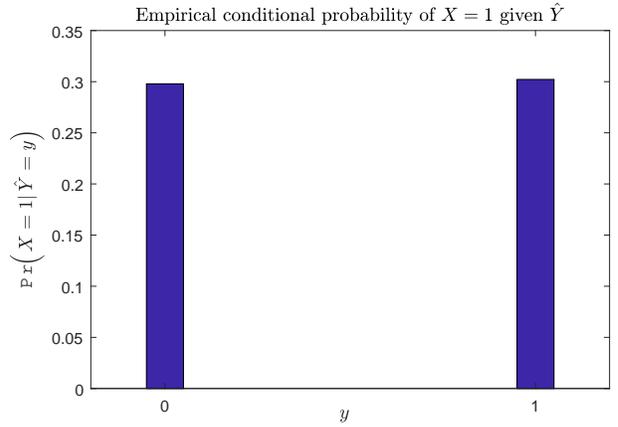}}}
\caption{The empirical conditional probability of $X=0$ ($a$) and $X=1$ ($b$) given $\hat{Y}=y$ under the optimal perfect privacy estimator of $Y$. The number of sensors is equal to $10$.}\label{Fig: F3}
\end{figure}

We next numerically study the dependency of the output of the optimal perfect privacy estimator on the private random variable. Fig. \ref{Fig: F3} shows the empirical conditional probability of $X$ given the output of the optimal perfect privacy estimator when the number of sensors is equal to $10$. According to this figure, under the perfect privacy estimator, the conditional probability of $X=1$ ($X=0$) does not depend on the estimated value of $Y$ and is equal to the probability of the event $X=1$ ($X=0$). This observation is in accordance with the fact that, under the perfect privacy constraint, the output of the estimator is independent of the private random variable.  Thus, an untrusted user is not able to infer the private information using the output of the estimator. 

\begin{figure}[!htp]
\centering{\includegraphics[scale=0.6]{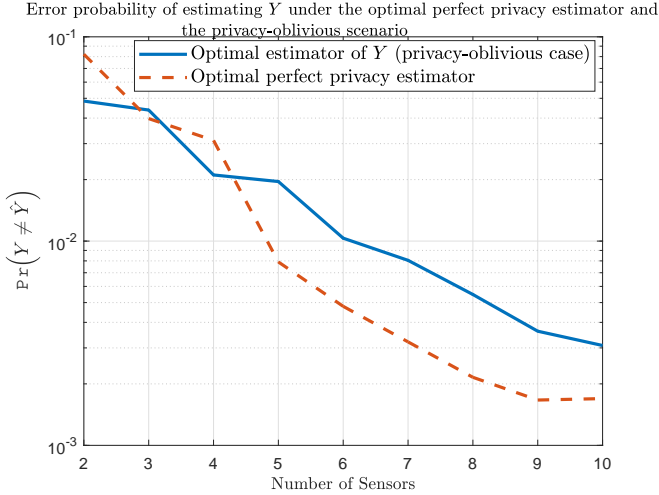}}
\caption{The error probability of estimating $Y$ under the optimal perfect privacy estimator and the optimal estimator of $Y$ in the privacy-oblivious scenario with the number of sensors.}\label{Fig: F4}
\end{figure}

Finally, in Fig. \ref{Fig: F4}, we compare the performance of the optimal perfect privacy estimator in recovering $Y$ with that of the privacy-oblivious estimation scenario. According to this figure, the performance of the optimal perfect privacy estimator in recovering the public random variable $Y$ is close to that of the optimal estimator of $Y$ in the privacy-oblivious scenario. Based on Fig. \ref{Fig: F2} and \ref{Fig: F3}, the perfect privacy estimator ensures that the  estimate of $Y$ does not contain any information about $X$ and its estimation reliability is close to the optimal estimator in the privacy-oblivious scenario, in the considered set-up.

\section{Conclusions}\label{Sec: Conc}
In this paper, we studied the optimal privacy-aware estimation of a public random variable when the sensors' measurements contain noisy information about the public random variable as well as a private random variable. The optimal estimation of the public random variable, under a binary loss function, with a constraint on the privacy level of the private random variable was studied. The conditional discrete entropy of the private random variable subject to the output of estimator was considered as the privacy metric and it was shown that the optimal estimator can be obtained by solving  (a possibly infinite dimensional) convex optimization problem. It was shown that the optimal perfect privacy estimator can be obtained by solving a convex optimization problem. 
\appendices
\section{Proof of Lemma \ref{Lem: Obj}}\label{App: Obj}
We first write the average loss function as 
\begin{align}
\ES{L\paren{Y,\hat{Y}_P\paren{Z}}}&=\ES{\I{Y\neq\hat{Y}_P\paren{Z}}}\nonumber\\
&=\PRP{Y\neq\hat{Y}_P\paren{Z}}\nonumber\\
&=1-\PRP{Y=\hat{Y}_P\paren{Z}}\nonumber
\end{align}
Using the Bayes' rule, $\PRP{Y=\hat{Y}_P\paren{Z}}$ can be written as \eqref{Eq: PE-Ex} where $p_Z\paren{z}$ is the probability density function of $Z$ and $(a)$ follows from the fact that, given $Z=z$, $\hat{Y}_P\paren{Z}$ is independent of $Y$. 
\begin{figure*}
\begin{align}\label{Eq: PE-Ex}
\PRP{Y=\hat{Y}_P\paren{Z}}&=\int \CP{Y=\hat{Y}_P\paren{Z}}{Z=z}p_Z\paren{z}dz\nonumber\\
&=\int \sum_i\CP{Y=\hat{Y}_P\paren{Z},Y=y_i}{Z=z}p_Z\paren{z}dz\nonumber\\
&=\int \sum_i\CP{Y=\hat{Y}_P\paren{Z}}{Z=z,Y=y_i}\CP{Y=y_i}{Z=z}p_Z\paren{z}dz\nonumber\\
&=\int \sum_i\CP{\hat{Y}_P\paren{Z}=y_i}{Z=z,Y=y_i}\CP{Y=y_i}{Z=z}p_Z\paren{z}dz\nonumber\\
&\stackrel{(a)}{=}\int \sum_i\CP{\hat{Y}_P\paren{Z}=y_i}{Z=z}\CP{Y=y_i}{Z=z}p_Z\paren{z}dz\nonumber\\
&=\sum_i\int P_i\paren{z}\CP{Y=y_i}{Z=z}p_Z\paren{z}dz
\end{align}
\hrule
\end{figure*}
\section{Proof of Lemma \ref{App: P-Conv}}\label{App: P-Conv}
Using the definition of the conditional entropy, the privacy constraint can be written as \eqref{Eq: E-Exp}. 
\begin{figure*}
\begin{align}\label{Eq: E-Exp}
\CDSE{X}{\hat{Y}_P\paren{Z}}&=\DSE{X}-\MI{X}{\hat{Y}_P\paren{Z}}\nonumber\\
&=\DSE{X}-\sum_j\PRP{X=x_j}\KLD{\CPMF{\hat{Y}_P}{y}{X=x_j}}{\PMF{\hat{Y}_P}{y}}
\end{align}
\hrule
\end{figure*}
 The probability mass functions $\PMF{\hat{Y}_P}{y}$ and  $\CPMF{\hat{Y}_P}{y}{X=x_j}$ can be expanded as 
\begin{align}
\PMF{\hat{Y}_P}{y_i}&=\PRP{\hat{Y}_P\paren{Z}=y_i}\nonumber\\
&=\int\CP{\hat{Y}_P\paren{Z}=y_i}{Z=z}p_Z\paren{z}dz\nonumber\\
&=\int P_i\paren{z}p_Z\paren{z}dz
\end{align}
and 
\begin{align}
&\CPMF{\hat{Y}_P}{y_i}{X=x_j}\nonumber\\
&=\CP{\hat{Y}_P\paren{Z}=y_i}{X=x_j}\nonumber\\
&=\int\CP{\hat{Y}_P\paren{Z}=y_i}{X=x_j,Z=z}\CD{Z}{z}{X=x_j}dz\nonumber\\
&=\int\CP{\hat{Y}_P\paren{Z}=y_i}{Z=z}\CD{Z}{z}{X=x_j}dz\nonumber\\
&=\int P_i\paren{z}\CD{Z}{z}{X=x_j}dz\nonumber
\end{align}
 respectively. Let $P^\prime\paren{z}=\left[P^\prime_i\paren{z}\right]_i$ denote a set of positive functions defined on $\R$ with $\sum_iP^\prime_i\paren{z}=1$ for all $z\in\R$. Consider a new estimator of $Y$, denoted by $\hat{Y}_{\tilde{P}}\paren{Z}$, where $\tilde{P}\paren{z}$ is a convex combination of $P\paren{z}$ and $P^\prime\paren{z}$, that is $\tilde{P}\paren{z}=\left[\alpha P_i\paren{z}+\paren{1-\alpha}P^\prime_i\paren{z}\right]_i$.  Thus, we have 
\begin{align}
\PMF{\hat{Y}_{\tilde{P}}}{y_i}&=\alpha \int P_i\paren{z}p_Z\paren{z}dz+\paren{1-\alpha} \int P^\prime_i\paren{z}p_Z\paren{z}dz\nonumber\\
&=\alpha \PMF{\hat{Y}_P}{y_i}+\paren{1-\alpha}\PMF{\hat{Y}_{P^\prime}}{y_i}
\end{align}
and, 
\begin{align}
\CPMF{\hat{Y}_{\tilde{P}}}{y_i}{X=x_j}&=\alpha \int P_i\paren{z}\CD{Z}{z}{X=x_j}dz\nonumber\\
&+\paren{1-\alpha} \int P^\prime_i\paren{z}\CD{Z}{z}{X=x_j}dz\nonumber\\
&=\alpha \CPMF{\hat{Y}_P}{y_i}{X=x_j}\nonumber\\
\hspace{1cm}&+\paren{1-\alpha} \CPMF{\hat{Y}_{P^\prime}}{y_i}{X=x_j}
\end{align}
where $\hat{Y}_{P^\prime}\paren{Z}$ is the estimator of $Y$ based on $P^\prime\paren{z}$. Using the convexity of KL divergence \cite{CT06}, $\KLD{\CPMF{\hat{Y}_{\tilde{P}}}{y}{X=x_j}}{\PMF{\hat{Y}_{\tilde{P}}}{y}}$ can be upper bounded as  \eqref{Eq: KL-conv}, for $j=\left\{1,\cdots,n\right\}$.
\begin{figure*}
\begin{align}\label{Eq: KL-conv}
\KLD{\CPMF{\hat{Y}_{\tilde{P}}}{y}{X=x_j}}{\PMF{\hat{Y}_{\tilde{P}}}{y}}&=\KLD{\alpha\CPMF{\hat{Y}_P}{y}{X=x_j}+\paren{1-\alpha}\CPMF{\hat{Y}_{P^\prime}}{y}{X=x_j}}{\alpha\PMF{\hat{Y}_P}{y}+\paren{1-\alpha}\PMF{\hat{Y}_{P^\prime}}{y}}\nonumber\\
&\leq \alpha \KLD{\CPMF{\hat{Y}_P}{y}{X=x_j}}{\alpha\PMF{\hat{Y}_P}{y}}+\paren{1-\alpha}\KLD{\CPMF{\hat{Y}_{P^\prime}}{y}{X=x_j}}{\PMF{\hat{Y}_{P^\prime}}{y}}
\end{align}
\hrule
\end{figure*}
which implies that the privacy constraint is convex in $P\paren{z}$.
\section{Proof of Lemma \ref{Lem: KKT}}\label{App: KKT}
Note that the mutual information between $X$ and $\hat{Y}_P\paren{Z}$ can be expanded as \eqref{Eq: MI-Exp}.
\begin{figure*}
\begin{align}\label{Eq: MI-Exp}
\MI{X}{\hat{Y}_P\paren{Z}}&=\sum_j\PRP{X=x_j}\KLD{\CPMF{\hat{Y}_P}{y}{X=x_j}}{\PMF{\hat{Y}_P}{y}}\nonumber\\
&=\sum_j\PRP{X=x_j}\sum_i\paren{\sum_lP_{il}\CP{Z\in B_l}{X=x_j}}\log\frac{\sum_lP_{il}\CP{Z\in B_l}{X=x_j}}{\sum_lP_{il}\PRP{Z\in B_l}}
\end{align}
\hrule
\end{figure*}
Using \eqref{Eq: MI-Exp}, the Lagrangian of the optimization problem \eqref{Eq: OP-D} can be written as \eqref{Eq: Lagran} where $\vec{\lambda}$ is the vector of Lagrange multipliers associated with the equality constraints, $\mu$ is the Lagrange multiplier associated with the privacy constraint and $\vec{P}=\left\{P_{ij}\right\}_{ij}$.
\begin{figure*}
\begin{align}\label{Eq: Lagran}
L\paren{\vec{P},\vec{\lambda},\mu}&=\sum_{i=1}^m\sum_{l=1}^NP_{il}\CP{Z\in B_l}{Y=y_i}\PRP{Y=y_i}\nonumber\\
&-\mu \sum_j\PRP{X=x_j}\sum_i\paren{\sum_lP_{il}\CP{Z\in B_l}{X=x_j}}\log\frac{\sum_lP_{il}\CP{Z\in B_l}{X=x_j}}{\sum_lP_{il}\PRP{Z\in B_l}}\nonumber\\
&+\sum_{l}\lambda_l\paren{\sum_iP_{il}-1}
\end{align}
\hrule
\end{figure*}

The partial derivative of Lagrangian with respect to $P_{hk}$ can be written as 
\begin{figure*}
\begin{align}
\frac{\partial L\paren{\vec{P},\vec{\lambda},\mu}}{\partial P_{hk}}&=\CP{Z\in B_k}{Y=y_h}\PRP{Y=y_h}+\lambda_k\nonumber\\
&-\mu\sum_j\PRP{X=x_j}\left[\CP{Z\in B_k}{X=x_j}\log\frac{\sum_l P_{hl}\CP{Z\in B_l}{X=x_j}}{\sum_{l}P_{hl}\PRP{Z\in B_l}}\right.\nonumber\\
&-\left.\paren{\sum_lP_{hl}\CP{Z \in B_l}{X=x_j}}\paren{\frac{\CP{Z\in B_k}{X=x_j}}{\sum_lP_{hl}\CP{Z\in B_l}{X=x_j}}-\frac{\PRP{Z\in B_k}}{\sum_l P_{hl}\PRP{Z \in B_l}}}\right] \nonumber\\
&\stackrel{(a)}{=}\CP{Z\in B_k}{Y=y_h}\PRP{Y=y_h}+\lambda_k\nonumber\\
&-\mu\sum_j\PRP{X=x_j}\left[\CP{Z\in B_k}{X=x_j}\log\frac{\sum_l P_{hl}\CP{Z\in B_l}{X=x_j}}{\sum_{l}P_{hl}\PRP{Z\in B_l}}\right]
\end{align}
\hrule
\end{figure*}
where $(a)$ follows from \eqref{Eq: Zero-Lag}. 
\begin{figure*}
\begin{align}\label{Eq: Zero-Lag}
\sum_j\PRP{X=x_j}\paren{\sum_lP_{hl}\CP{Z \in B_l}{X=x_j}}\paren{\frac{\CP{Z\in B_k}{X=x_j}}{\sum_lP_{hl}\CP{Z\in B_l}{X=x_j}}-\frac{\PRP{Z\in B_k}}{\sum_l P_{hl}\PRP{Z \in B_l}}}=0
\end{align}
\hrule
\end{figure*}

Note that the optimization problem  \eqref{Eq: OP-D} is a convex optimization problem and it is straightforward to show that the Slater's condition holds for this problem. Thus, using the necessary and sufficient Karush-Kuhn-Tucker (KKT) conditions, we have 
\begin{align}
\frac{\partial L\paren{\vec{P},\vec{\lambda},\mu}}{\partial P_{hk}}|_{\vec{P}^\star,\vec{\lambda^\star,\mu^\star}}\left\{
\begin{array}{cc}
=0 \quad \text{if}&  0<P^\star_{hk}<1\nonumber\\
\leq  0  \quad\text{if}& P^\star_{hk}=0\nonumber\\
\geq 0 \quad \text{if}& P^\star_{hk}=1\nonumber
\end{array}
\right.
\end{align}
where $\vec{P}^\star$ is the optimal solution, and $\vec{\lambda}^\star$ and $\mu^\star$ are the dual optimal variables.
\section{Proof of Lemma \ref{Lem: PP-C-CM}}\label{App: PP-C-DM}
To satisfy the perfect privacy condition, we need to have  
\begin{align}\label{Eq: Ind}
\PRP{\hat{Y}\paren{Z}=y_i, X=x_j}=\PRP{\hat{Y}\paren{Z}=y_i}\PRP{X=x_j} \forall i,j\nonumber
\end{align}
Note that $\PRP{\hat{Y}\paren{Z}=y_i,X=x_j}$ and $\PRP{\hat{Y}\paren{Z}=y_i}$ can be expanded as
\begin{align}
&\PRP{X=x_j,\hat{Y}\paren{Z}=y_i}\nonumber\\
&=\sum_l{\PRP{X=x_j,\hat{Y}\paren{Z}=y_i,Z\in B_l}}\nonumber\\
&=\PRP{X=x_j}\sum_l\CP{\hat{Y}\paren{Z}=y_i}{Z\in B_l}\CP{Z\in B_l}{X=x_j}\nonumber\\
&=\PRP{X=x_j}\sum_lP_{il}\CP{Z\in B_l}{X=x_j}\nonumber
\end{align}
and 
\begin{align}
&\PRP{\hat{Y}\paren{Z}=y_i}\nonumber\\
&=\sum_l\PRP{\hat{Y}\paren{Z}=y_i,Z\in B_l}\nonumber\\
&=\sum_l\CP{\hat{Y}\paren{Z}=y_i}{Z\in B_l}\PRP{Z \in B_l}\nonumber\\
&=\sum_lP_{il}\PRP{Z \in B_l}\nonumber
\end{align}
Thus, for $\PRP{X=x_j}\neq 0$, the perfect privacy requirement  can be expressed as 
\begin{align}
\sum_l P_{il}\paren{\CP{Z\in B_l}{X=x_j}-\PRP{Z \in B_l}}=0 \quad \forall i, j\nonumber
\end{align}
For a given $i$, the condition above can be expressed as 
\begin{align}
\Phi\vec{P}_i=\vec{0}
\end{align}
which implies that $\vec{P}_i\in \Null{\Phi}$.
\section{Proof of Lemma \ref{Lem: PP-Feas}}\label{App: PP-DM-Feas}
 Let $\mathcal{P}=\mathcal{P}_1\times\cdots\times\mathcal{P}_N$ denote the joint probability simplex corresponding to the randomization probabilities of bins, \emph{i.e.,} $\paren{P_{1l},\cdots,P_{ml}}^\top\in\mathcal{P}_l$ for all $l\in\left\{1,\cdots,N\right\}$. We show that there are infinity many points in $\mathcal{P}$ which satisfy the perfect privacy condition if $N> n$.
 
 Note that for $N>n$, the null space of $\Phi$ is non-empty. Pick a set of positive real numbers $\left\{\lambda_i \right\}_{i=1}^m$  where $0< \lambda_i<1$ for all $i$ and  $\sum_{i}\lambda_i=1$. Let $\paren{\lambda_1,\cdots,\lambda_m}^\top$  be the randomization probability of bin $l$ for all $1\leq l\leq N$, \emph{i.e.,} $P_{il}=\lambda_i$ for all $i,l$. Thus, we have $\vec{P}_i=\paren{\lambda_i,\cdots,\lambda_i}$. Note that $\vec{P}_i$ belongs to $\Null{\Phi}$ as the sum of the elements in each row of $\Phi$ is equal to zero.
 \section{Proof of Lemma \ref{Lem: PP-CM-Feas}}\label{App: PP-CM-Feas}
 Pick $N>n$ and $\left\{B_i\subset \R\right\}_{i=1}^N$ where $B_i$s are arbitrary discretization bins. According to Lemma \ref{Lem: PP-Feas}, the feasible set of the optimization problem \eqref{Eq: OP-PP} is non-empty. Let  $\left\{\bar{P}_l\right\}_l$ denote a point in the feasible set of \eqref{Eq: OP-PP} where $\bar{P}_l=\paren{P_{1l},\cdots,P_{ml}}^\top$ denotes the randomization probabilities corresponding to the bin $l$. Now, we construct $\left\{P_i\paren{z}\right\}_i$ as follows. For $z\in B_l$, let $P_i(z)=P_{il}$ for all $i$. Thus, any feasible solution of \eqref{Eq: OP-PP} corresponds to a piecewise constant solution for \eqref{Eq: OP-PP-CM} which satisfies the privacy constraint.
\bibliographystyle{IEEEtran}
\bibliography{MTNS}

\end{document}